\documentclass[a4paper,12pt]{article}
\usepackage{amsmath, amssymb, amsthm}
\usepackage[dvips,final]{epsfig}
\usepackage[dvips,final]{graphicx}
\usepackage{amsmath}
\usepackage{amsthm}
\usepackage{enumerate}
\usepackage{setspace}
\usepackage{amsfonts}
\usepackage{empheq}
\usepackage{amssymb}

\usepackage{graphicx}
\usepackage[english]{babel}
\usepackage{color}
\definecolor{lightgray}{rgb}{0.83, 0.83, 0.83}
\definecolor{lightblue}{rgb}{0.83, 0.83, 1}
\definecolor{ao}{rgb}{0.0, 0.5, 0.0}
\definecolor{dg}{rgb}{0.01, 0.75, 0.24}
\definecolor{orange}{rgb}{1.0, 0.49, 0.0}
\definecolor{cardinal}{rgb}{0.77, 0.12, 0.23}
\definecolor{dollarbill}{rgb}{0.52, 0.73, 0.4}

\theoremstyle{plain} 
\newtheorem{thm}{Theorem}[section]

\theoremstyle{definition}

\theoremstyle{definition} 
\newtheorem{defn}{Definition}[section]

\theoremstyle{remark} 
\newtheorem{rmk}{Remark}

\newenvironment{michelarev}{\color{blue}}{\color{black}}
\newcommand{\bmicr}{\begin{michelarev}}
\newcommand{\emicr}{\end{michelarev}}

\newcommand{\fhi}{\varphi}

\def\dd{\,\mathrm{d}}
\def\dive{\mathrm{\,div\,}}
\def\ddt{\dfrac{d}{dt}}

\def\real{\mathbb{R}}

\def\AA{\mathcal{A}}
\def\SS{\mathcal{S}}
\def\PP{\mathcal{P}}
\def\CC{\mathcal{C}}
\def\BB{\mathcal{B}}

\def\io{\int_{\Omega}}

\def\itt{\int_0^{T}}

\def\12{\dfrac{1}{2}}
\def\p{\partial}

\numberwithin{equation}{section}

\begin{document}

\title{On a non-isothermal Cahn-Hilliard model \\ for tumor growth}
\author {Erica Ipocoana\thanks{University of Modena and Reggio Emilia, Dipartimento di Scienze Fisiche, Informatiche e Matematiche, via Campi 213/b, I-41125 Modena (Italy).\newline
E-mail: \textit{erica.ipocoana@unipr.it}}}

\maketitle
\maketitle
\begin{abstract}
We introduce here a new diffuse interface thermodynamically consistent non-isothermal model for tumor growth in presence of a nutrient in a domain $\Omega \subset \real^3$. In particular our system describes the growth of a tumor surrounded by healthy tissues, taking into account changes of temperature, proliferation of cells, nutrient consumption and apoptosis. Our aim consists in proving an existence result for our problem associated to the entropy formulation.
\end{abstract}


\section{Introduction}
The study of tumor growth processes has become of great interest also for mathematicians in recent years  \cite{am, bell, gatti, cl, ltz, ohp, wu}. Indeed, mathematical models might be able to give further insights in tumor growth behaviour. In particular, the framework of diffuse interface modeling with Cahn-Hilliard equations \cite{ch} has received increasing attention. In this context, the tumor is seen as an expanding mass surrounded by healthy tissues. Its evolution is assumed to be governed by mechanisms such as proliferation of cells via nutrient consumption, apoptosis \cite{fgr, jw, mrs19} and, in more complex models like \cite{garkewp, garken, garke, chem}, also chemotaxis and active transport of specific chemical species effects.
Moreover it is possible to include the effects of fluid flow into the evolution of the tumor, which brings to the so-called Cahn-Hilliard-Darcy models (see \cite{garke, jw}).
However, up to our knowledge it seems that even if the effects of variations of temperature have been studied for Cahn-Hilliard equations \cite{Eleuteri.Gatti.Schimperna2017, m.eleuterie.roccag.schimperna2015, iz}, they have been neglected so far in the analysis of tumor growth. From the medical point of view, the effects of temperature on tumor growth have not been completely understood yet, although they have been investigated since the very beginning of the 20th century \cite{5}. The general tendency of the scientific community seems to support the thesis that hyperthermia can lead to partial or complete destruction of tumor cells \cite{1,3,4,6}. In fact, it has also been observed that low ambient temperature influences the production of particular nutrients for the tumor \cite{2}. Nevertheless, we focus here on the case which does not take into account the production of a nutrient due to temperature.
In this work we introduce a new diffuse interface model for tumor growth, taking into account proliferation of cells, nutrient consumption and apoptosis and moreover temperature effects. Our aim consists in proving an existence result for weak entropy solutions (cnfr. Definition \ref{defensol}) to our model.
We remark that a rigorous mathematical theory of well-posedness results has been addressed in multiple works, such as \cite{fgr, garkewp, mrs19}.
From the biological point of view, we assume that tumor cells only die by apoptosis, therefore we do not take into account the possibility of tumor necrosis (differently e.g. from \cite{garken}). We also suppose that the healthy cells surrounding the tumor do not interact with the tumor itself, neglecting the possible response of the immune system.\\
According to these considerations, we will derive the following PDE system, describing the behaviour of a two-component mixture consisting of healthy cells and tumor cells
\begin{align}\label{eqfhi-i}
&\fhi_t = \Delta \mu +(\PP\sigma-\AA)h(\fhi)\\\label{eqmu-i}
&\mu = -\varepsilon\Delta\fhi + \dfrac{1}{\varepsilon}F'(\fhi)-\theta -\chi_\fhi \sigma\\\label{eqtemp-i}
&\theta_t +\theta\fhi_t -\dive(\kappa(\theta)\nabla\theta) = |\nabla\mu|^2\\\label{eqnutr-i}
&\sigma_t = \Delta \sigma -\CC\sigma h(\fhi)+\BB(\sigma_B-\sigma).
\end{align}
We carry out our analysis in $\Omega \times (0,\infty)$, where $\Omega \subset \real^3$ is a smooth domain. According to the derivation of the model shown in Section \ref{deriv}, we suppose that the system is isolated from the exterior. This condition translates in no-flux boundary conditions (i.e. homogeneous Neumann) for all the unknowns.\\
The evolution of the tumor is described by the order parameter $\fhi$  which represents the local concentration of tumor cells, $\fhi \in [-1,1]$, with $\lbrace\fhi = 1 \rbrace$ representing the tumor phase and $\lbrace\fhi = -1\rbrace$ the healthy one. Moreover $\mu$ denotes the chemical potential of phase transition from healthy to tumor cells, $\theta$ is the absolute temperature, $\kappa(\theta)$  represents the heat conductivity and $\varepsilon$ is a small parameter related to the thickness of interfacial layers. We denote by $\sigma$ the concentration of a nutrient consumed (only) by the tumor cells  (e.g. oxygen and glucose). The parameter $\chi_\fhi \geq 0$ is linked to transport mechanisms such as chemotaxis and active uptake. Althought we will show in Section \ref{deriv} how this parameter is included in the model, for sake of simplicity we will neglect it throughout the mathematical analysis, with the aim of including it in future works. The positive constant parameters $\PP, \AA, \CC$ and $\BB$ indicate respectively the tumor proliferation rate, apoptosis rate, nutrient consuption rate and nutrient supply rate.
The function $h$ is chosen as monotone increasing, nonnegative in $[-1,1]$ and such that $h(-1) = 0$ and $ h(1) = 1$. The tumor growth is thus described by the term $\PP \sigma h(\fhi)$, which reasonably increases proportionally to the concentration of tumor cells, while the death of tumor cells is modelled by the term $\AA h(\fhi)$. 
Therefore, according to \eqref{eqfhi-i}, if $\PP\sigma -\AA >0$, then the tumor expands and it happens faster when the concentration of tumor cells is already high. If otherwise $\PP\sigma -\AA <0$ then the tumor reduces and the tumor cells die faster when the concentration of tumor cells is large.
The term $\CC \sigma h(\fhi)$ represents the consumption of the nutrient by the tumor cells.
The term $\BB(\sigma-\sigma_B)$ is due to the fact that we consider here the case where the tumor has its own vasculature (as in e.g. \cite{nutu}, \cite{mrs19}), where the threshold $\sigma_B \in (0,1)$ is the constant nutrient concentration in the pre-existing vasculature. In particular, if $\sigma_B>\sigma$, $\BB(\sigma_B-\sigma)$ models the supply of nutrient from the blood vessels, on the other hand if $\sigma_B<\sigma$, $\BB(\sigma_B-\sigma)$ represents the transport of nutrient away from the domain.
Eventually, the function $F(s)$ represents a polynomial potential having at least cubic growth at infinity, whose assumptions will be specified in Section \ref{ass}. A simple choice might be a \textit{double-well potential} with equal minima at $s=\pm 1$ penalizing the deviation of the length $|\fhi|$ from its natural value $1$. This more general potential allows $\fhi$ to take values also outside of the significance interval $[-1,1]$, therefore we will carry out our analysis also in the case $|\fhi|>1$ and correspondingly extend function $h$. We also remark that although among Cahn-Hilliard literature the singular potentials, such as logarithm type (see e.g. \cite{fg}), are very common, the growth conditions that the problem requires make them unsuitable for our case, as it will be clear in Section \ref{ape}.\\
In this work we derive a new phase field model according to the laws of thermodynamics describing the tumor growth. The novelty of this contribution is to include possible variations of temperature in the model. The presence of nutrient concentration $\sigma$ in the system implies that here the spatial mean of $\fhi$ is not conserved in time (as we can see from equation \eqref{eqfhi-i}), therefore the derivation of the model cannot follow the standard techniques proposed e.g. in \cite{m.eleuterie.roccag.schimperna2015}. However we are able to gain enough regularity for the quadruple $(\fhi, \mu,\theta,\sigma)$ in order to prove the existence of weak solutions to the initial-boundary value problem associated to \eqref{eqfhi-i}--\eqref{eqnutr-i}.\\
The structure of this paper is the following. In Section \ref{deriv} we derive system \eqref{eqfhi-i}--\eqref{eqnutr-i} according to the approach proposed by Gurtin in \cite{gurtin}. Then we proceed with the mathematical analysis of our problem in the case $\varepsilon =1, \;\chi_\fhi =0$. In particular, Section \ref{exreg} is devoted to give the setting and to present the main result of this work (which is Theorem \ref{maint}) concerning the existence of weak entropy solutions to our problem. The proof is carried out in two steps. In Section \ref{ape} we gain a priori bounds for $(\fhi, \mu,\theta,\sigma)$. In Section \ref{wss} we use the weak sequential stability argument to prove the existence of weak entropy solutions. Namely, we exploit the a priori bounds obtained for a sequence of weak entropy solutions together with standard compactness results to pass to the limit. 

\section{Derivation of the model}
\label{deriv}
We suppose that a two-component mixture consisting of healthy cells and tumor cells occupies an open spatial domain $\Omega \subset \real^3$. We denote by $\fhi(x,t)$ the tumor phase concentration, $\theta(x,t)$ is the absolute temperature and $\sigma(x,t)$ is the concentration of a nutrient for the tumor cells.
According  to the Ginzburg-Landau theory for phase transitions, we postulate the free energy density $\psi$ in the form
\begin{eqnarray}\label{fren}
\psi = \dfrac{\varepsilon}{2}|\nabla\fhi|^2 + \dfrac{1}{\varepsilon}F(\fhi)+f(\theta)-\theta \fhi + N(\fhi,\sigma).
\end{eqnarray} 
Here, $\varepsilon$ is a positive constant depending on the interface thickness. The function $F$ represents a polynomial potential having at least cubic growth at infinity. The easiest choice is taking $F(\fhi) = (|\fhi|^2-1)^2$, known in the literature as the \textit{double-well potential}.\\
The term $f$ in \eqref{fren} describes the part of free energy which is purely caloric and is related to the specific heat $c_V(\theta) = Q'(\theta)$ through relation $Q(\theta)=f(\theta)-\theta f'(\theta)$. In the following we assume the specific heat $c_V\equiv 1$. Moreover we recall that it holds
\begin{eqnarray} \label{q}
q = Q\theta
\end{eqnarray}
where $q$ denotes the heat flux.
Eventually, the latter term $N$ in equation \eqref{fren} describes both the chemical energy of the nutrient and the energy contributions given by the interactions between the tumor tissues and the nutrient.\\ 
One of the main difficulties we have to afford in the derivation of our model is that, differently from standard Cahn-Hilliard models (such as \cite{m.eleuterie.roccag.schimperna2015}), the spatial mean of the tumor phase concentration $\fhi$ is not conserved. Indeed the tumor may grow or shrink according to the right hand side of \eqref{eqfhi-i}. In order to deal with this issue, we follow Gurtin's approach (used e.g. in \cite{ms, m02}) proposed in \cite{gurtin}, namely we treat separately the balance laws and the constitutive relations, moreover we introduce the following new balance law for internal microforces
\begin{eqnarray}\label{micro}
\dive{\zeta}+\pi = 0,
\end{eqnarray}
where $\zeta$ is a vector representing the microstress and $\pi$ is a scalar corresponding to the internal microforces. We remark that here we neglect the external actions.\\
The mass balance law reads
\begin{eqnarray}\label{mass}
\fhi_t = -\dive h + m,
\end{eqnarray}
where $h$ is the mass flux and $m$ is the external mass supply. Moreover the internal energy density of the system is given by
\begin{eqnarray}\label{inten}
e = \psi + \theta s.
\end{eqnarray}
Here, $s$ denotes the entropy of the system, which has the following expression, according to \eqref{fren}
\begin{eqnarray}\label{entropy}
s = -\dfrac{\partial\psi}{\partial\theta} = -f'(\theta)+\fhi.
\end{eqnarray}
Combining the two previous formulas we infer
\begin{eqnarray}\label{eder}
\dfrac{\p e}{\p t} = \dfrac{\p \psi}{\p t}+\theta \dfrac{\p s}{\p t}+s\dfrac{\p \theta}{\p t} 
= \dfrac{\p \psi}{\p \fhi}\dfrac{\p \fhi}{\p t}+\dfrac{\p \psi}{\p \nabla \fhi}\dfrac{\p \nabla \fhi}{\p t} +\theta \dfrac{\p s}{\p t}
\end{eqnarray}
and consequently
\begin{eqnarray}\label{etrick}
\dfrac{\p \psi}{\p t}+s\dfrac{\p \theta}{\p t} 
= \dfrac{\p \psi}{\p \fhi}\dfrac{\p \fhi}{\p t}+\dfrac{\p \psi}{\p \nabla \fhi}\dfrac{\p \nabla \fhi}{\p t}.
\end{eqnarray}
\medskip
\noindent \textbf{Cahn-Hilliard system}\\
The derivation of our system is based on the two fundamental laws of thermodynamics.
According to \cite{gurtin}, we write the first law in the form
\begin{eqnarray}\label{1tif}
\dfrac{d}{dt}\int_{\textit{R}} e \dd x = -\int_{\partial \textit{R}} q\cdot \nu \dd \eta + \mathcal{W}(\textit{R})+\mathcal{M}(\textit{R}),
\end{eqnarray}
where \textit{R} is the control volume, $\nu$ is the outward unit normal to $\partial R$ and
\begin{align*}
&\mathcal{W}(\textit{R}) = \int_{\partial \textit{R}} (\zeta \cdot \nu)\dfrac{\partial \fhi}{\partial t} \dd \eta,\\
&\mathcal{M}(\textit{R}) = -\int_{\partial \textit{R}} \mu h\cdot\nu\dd \eta + \int_{\textit{R}}\mu m \dd x
\end{align*}
are the rate of working and the rate at which free energy is added to \textit{R} (assuming no heat supply) respectively.
Using Green's formula, we can rewrite \eqref{1tif} as
\begin{eqnarray}\label{1tdf}
\dfrac{\partial e}{\partial t}= -\dive q + \dfrac{\partial\fhi}{\partial t} \dive \zeta + \zeta\cdot \nabla\dfrac{\partial \fhi}{\partial t} - h \cdot\nabla \mu-\mu \dive h +\mu m.
\end{eqnarray}
Since the control volume \textit{R} is arbitary, exploiting the mass balance \eqref{mass} and the microforce balance \eqref{micro}, we infer
\begin{eqnarray} \label{2tdf}
\dfrac{\partial e}{\partial t}= -\dive q + (\mu-\pi)\dfrac{\partial\fhi}{\partial t}+ \zeta\cdot \nabla\dfrac{\partial \fhi}{\partial t} - h \nabla \mu.
\end{eqnarray}
We now impose the validity of the second law of thermodinamics in the form of the Clausius-Duhem inequality
\begin{eqnarray}\label{C-D}
\theta \left( \dfrac{\partial s}{\partial t} + \dive Q\right) \geq 0.
\end{eqnarray}
We develop the left hand side of \eqref{C-D} as follows
\begin{small}
\begin{align*}
\theta \left( \dfrac{\partial s}{\partial t} + \dive Q\right)
& \stackrel{\eqref{inten}}{=}\dfrac{\partial e }{\partial t}-\dfrac{\partial \psi}{\partial t}-s\dfrac{\partial \theta}{\partial t} +\theta\dive Q\\
&\stackrel{\eqref{q}}{=}\dfrac{\p e}{\p t} -\dfrac{\p \psi}{\p t}-s\dfrac{\p \theta}{\p t} + \dive q -Q\cdot\nabla \theta\\
& \stackrel{\eqref{2tdf}}{=}(\mu-\pi)\dfrac{\p \fhi}{\p t}+\zeta \cdot\nabla \dfrac{\p\fhi}{\p t}-h\nabla\mu -\dfrac{\p\psi}{\p t}-s\dfrac{\p \theta}{\p t}-Q\cdot\nabla\theta\\
&\stackrel{\eqref{etrick}}{=}\left(\mu-\pi-\dfrac{\p\psi}{\p\fhi}\right)\dfrac{\p \fhi}{\p t}+\left(\zeta-\dfrac{\p \psi}{\p \nabla\fhi}\right)\dfrac{\p\nabla\fhi}{\p t} -h\nabla\mu -Q\cdot\nabla\theta.
\end{align*}
\end{small}
In order to satisfy relation \eqref{C-D}, we impose
\begin{align}\label{mu}
&\mu-\pi-\dfrac{\p\psi}{\p\fhi}=0,\\ \label{zeta}
&\zeta =\dfrac{\p \psi}{\p\nabla\fhi},\\ \label{gradmu}
& h\nabla\mu +Q\cdot\nabla\theta \leq 0,
\end{align}
where in particular in order for \eqref{gradmu} to hold, we exploited Fourier's law 
\begin{equation}\label{fourier}
q = -\kappa(\theta)\nabla\theta, 
\end{equation}
with $\kappa = \kappa(\theta)>0$ heat conductivity.\\
The combination of \eqref{fren} and \eqref{zeta} straightly gives
\begin{eqnarray}\label{zetaespl}
\zeta = \varepsilon\nabla\fhi
\end{eqnarray} 
which leads to, according to \eqref{fren}, \eqref{micro} and \eqref{mu},
\begin{eqnarray} \label{muespl} 
\mu = -\varepsilon\Delta\fhi + \dfrac{1}{\varepsilon}F'(\fhi)-\theta +\dfrac{\p N}{\p \fhi}.
\end{eqnarray}
Eventually, inequality \eqref{gradmu} can be satisfied choosing $h = -\nabla\mu$, which is a suitable assumption according to \cite{gurtin}. Therefore equation \eqref{mass} reads
\begin{eqnarray} \label{fhi}
\fhi_t = \Delta \mu +m.
\end{eqnarray}

\medskip
\noindent\textbf{Temperature equation.}\\
We start from the internal energy equation \eqref{2tdf}, taking advantage of \eqref{zeta} and of the expression for the chemical potential \eqref{mu}, therefore
\begin{align*}
\dfrac{\p e}{ \p t} = -\dive q +\dfrac{\p \psi}{\p \fhi}\dfrac{\p\fhi}{\p t}+\dfrac{\p\psi}{\p\nabla\fhi}\dfrac{\p\nabla\fhi}{\p t}-h\nabla\mu.
\end{align*}
Now, exploiting the assumption $h =-\nabla \mu$ and Fourier's law \eqref{fourier}, we infer
\begin{align*}
\dfrac{\p e}{ \p t} = \dive(\kappa(\theta)\nabla\theta)+\dfrac{\p \psi}{\p \fhi}\dfrac{\p\fhi}{\p t}+\dfrac{\p\psi}{\p\nabla\fhi}\dfrac{\p\nabla\fhi}{\p t}+|\nabla\mu|^2
\end{align*}
and by identity \eqref{eder},
\begin{align}\label{fstemp}
\theta\dfrac{\p s}{\p t} -\dive(\kappa(\theta)\nabla\theta) = |\nabla\mu|^2.
\end{align}
From \eqref{entropy}, we might write 
\[\theta \frac{\partial s}{\partial t} = \theta(-f'(\theta))_t+\theta \fhi_t.\]
On the other hand, according to the definition of $Q$, it holds $(Q(\theta))_t =(-f'(\theta))_t$, with in particular $(Q(\theta))_t = Q'(\theta)\theta_t$. Since we supposed that we are considering the case in which the specific heat $c_V=1$,
it follows that $Q'(\theta) =1$. This implies that
\[\theta s_t = \theta_t +\theta \fhi_t.\] 
Thus, equation \eqref{fstemp} reads
\begin{eqnarray} \label{tempeq}
\theta_t +\theta\fhi_t -\dive(\kappa(\theta)\nabla\theta) = |\nabla\mu|^2.
\end{eqnarray}
\medskip\\
\textbf{Nutrient equation.}\\
We postulate the nutrient balance equation in the form
\begin{eqnarray} \label{nutrbal}
\sigma_t = -\dive J-\mathcal{S},
\end{eqnarray}
where $J$ is the nutrient flux and $\mathcal{S}$ denotes a source/sink term for the nutrient. Motivated by \cite{garke}, we choose $J = -\nabla\sigma$, therefore equation \eqref{nutrbal} reads
\begin{eqnarray} \label{nutr}
\sigma_t = \Delta \sigma -\mathcal{S}.
\end{eqnarray}

\subsection{Constitutive relations}
Owing to \cite{nutu, garke, mrs19}, we now make the following constitutive assumptions.
\begin{itemize} 
\item[$\bullet$] $m = (\PP\sigma-\AA)h(\fhi)$,\\ 
where $h(\fhi)$ is a monotone increasing, nonnegative function in $[-1,1]$ and such that $h(-1) = 0$ and $ h(1) = 1$. Hence this relation states that on one hand the tumor growth is proportional to the nutrient supply in the tumoral region. This assumption reflects the fact that it often happens that tumors bring mutations which switch off certain growth inhibiting proteins. Therefore the tumor cells increasing is limited only by the supply of nutrients, despite of healthy cells where the mitotic cycle regulates the growth. On the other hand, when we are in the healthy region, this equation shows that the proliferation rate of the tumor is greater than the one of healthy cells.
\item[$\bullet$] $\dfrac{\p N}{\p \fhi} = -\chi_{\fhi}\sigma$, in fact, we take $\chi_{\fhi} =0$.\\ 
Indeed, this equation is due to the mechanism of chemotaxis, which we exclude in our analysis.
\item[$\bullet$] $\SS=\CC\sigma h(\fhi)-\BB(\sigma_B-\sigma)$.\\
We here assume that the sink/source of nutrient is regulated by consumption of nutrients and the term $\BB(\sigma_B-\sigma)$ which models the fact that we here consider the case in which the tumor has its own vasculature. In particular the threshold $\sigma_B$ indicates whether the nutrient is supplied to the tumor or transported away.
\end{itemize}

\section{Existence of solutions}
\label{exreg}
In this section we present the main result of this work, concerning the existence of solutions for the tumor growth model \eqref{eqfhi-i}--\eqref{eqnutr-i} for $\chi_\fhi =0$ and $\varepsilon = 1$. Namely, we work on system
\begin{align}\label{eqfhi}
&\fhi_t = \Delta \mu +(\PP\sigma-\AA)h(\fhi)\\\label{eqmu}
&\mu = -\Delta\fhi + F'(\fhi)-\theta \\\label{eqtemp}
&\theta_t +\theta\fhi_t -\dive(\kappa(\theta)\nabla\theta) = |\nabla\mu|^2\\\label{eqnutr}
&\sigma_t = \Delta \sigma -\CC\sigma h(\fhi)+\BB(\sigma_B-\sigma).
\end{align}

\subsection{Notation}
In order to carry out a mathematical analysis of our problem, let us introduce some notation we will use in the sequel.
We recall that $\Omega$ is a smooth domain of $\real^3$ and we denote by $\Gamma$ its boundary. For sake of semplicity, let us assume $|\Omega| =1$. We denote by $(0,T)$ an assigned but otherwise arbitrary time interval. We set $H:=L^2(\Omega)$ and $V:=H^1(\Omega)$ and we will use these symbols also referring to vector valued functions. The symbol $(\cdot,\cdot)$ will
indicate the standard scalar product in $H$, while $\langle \cdot,\cdot \rangle$ will stand for the duality between $V'$ and $V$.
We denote by $\|\cdot\|_X$ the norm in the generic Banach space $X$. For brevity we will write $\|\cdot\|$ instead of $\|\cdot\|_H$. Still for brevity, we omit the variables of integration. We
will specify them when there could be a misinterpretation.\\
For any function $v \in V$, we define
\begin{align} \label{mean}
v_\Omega:= \dfrac{1}{|\Omega|} \io v = \io v,
\end{align}
where the last equality holds since we assumed $|\Omega| = 1$.\\
We recall the Poincar\'{e}-Wirtinger inequality
\begin{align}\label{poin}
\|v-v_{\Omega}\|\leq c_{\Omega}\|\nabla v\| \quad \forall v\in V
\end{align}
and the non-linear Poincar\'{e} inequality
\begin{align}\label{nlp}
\|v^{\frac{p}{2}}\|^2_V\leq c_p\left(\|v\|^p_{L^1(\Omega)}+\| \nabla v ^{\frac{p}{2}}\|^2\right), 
\end{align}
which holds $\forall v \in L^1(\Omega)$ s.t. $\nabla v^{\frac{p}{2}}\in L^2(\Omega)$ and $\forall p\in [2,\infty)$.

\subsection{Assumptions}
\label{ass}
We assume the coefficients $\PP,\AA,\BB$ and $\CC$ to be strictly positive and $\sigma_B \in (0,1)$.
Next, we suppose that the derivative of potential $F \in C^1_{loc}(\real,\real)$ decomposes as a sum of a monotone increasing part $\beta$ and a linear perturbation, namely
\begin{align}\label{F}
F'(r)=\beta(r)-\lambda r \quad \lambda \geq 0, \; r\in\real.
\end{align}
Moreover we normalize $\beta$ s.t. $\beta(0)=0$ and we require
\begin{align} \label{beta}
\exists c_\beta >0 \; \;\textrm{s.t.} \;\; &|\beta(r)|\leq c_\beta(1+F(r)) \quad \forall r \in \real,\\\label{bsup}
&|\beta(r)|\geq k|r| \textrm{ for some } k>0,
\end{align}
where \eqref{beta} means that $F$ has at most an exponential growth at infinity, while \eqref{bsup} states that $\beta$ has superlinear growth.
Moreover, we assume potential $F$ to be strictly positive.\\
Next, we assume $h \in C^1(\real)$ increasingly monotone s.t. 
\begin{itemize}
\item[i)]$h(-1) = 0, \; h(r) \equiv 1 \quad \forall r\geq 1$. 
\item[ii)]$\exists \ \underline{h}\geq 0 \textrm{ and } \underline{\fhi}\leq -1 \textrm{  s.t. } h(r) \equiv -\underline{h} \quad \forall r\leq \underline{\fhi}$.
\end{itemize}
Therefore $h$ is globally Lipschitz continuous and there exists a constant $c>0$ s.t.
\begin{align}\label{hh'}
|h(r)|+|h'(r)|\leq c \quad \forall r \in \real.
\end{align}
Moreover we assume the thermal conductivity to depend on the absolute temperature $\theta$ as follows
\begin{align} \label{heatassumption}
\kappa(\theta) = 1 + \theta^q, \quad q \in [2, \infty), \quad \theta \geq 0.
\end{align} 
Eventually, we require the initial data to be such that
\begin{align}\nonumber
&\fhi|_{t=0} = \fhi_0, \quad \fhi_0 \in V, \quad F(\fhi_0) \in L^1(\Omega)
\\ \nonumber
& \theta|_{t=0} = \theta_0, \quad \; \theta_0 \in L^1(\Omega),\; \theta_0 >0 \textrm{ a.e. in }\Omega, \; \log \theta_0 \in L^1(\Omega) \\  \label{ic}
& \sigma|_{t=0} = \sigma_0, \quad \sigma_0 \in L^\infty(\Omega), \; 0\leq \sigma_0 \leq 1 \textrm{ a.e. in }\Omega 
\end{align}
where the last assumption on $\sigma_0$ is due to the interpretation of $\sigma$ as a nutrient concentration. We also recall that we couple our system with homogeneus Neumann boundary conditions for all the unknowns.

\subsection{Main result}

We here present what will be called a \textit{weak entropy solution}, already used e.g. in \cite{ms}, which is in fact weaker
than other corresponding notions appearing in related contexts. This is due to the fact that we do not get enough regularity to pass to the limit in some non-linear terms in the temperature equation \eqref{eqtemp}.\\
Multiplying \eqref{eqtemp} by $\frac{1}{\theta}$, we have
\begin{equation}\label{Leq}
(\Lambda(\theta) + \fhi)_t -\dive \left(\dfrac{\kappa (\theta)\nabla\theta}{\theta}\right)= \dfrac{\kappa(\theta)}{\theta^2}|
\nabla \theta|^2 + \dfrac{|\nabla\mu|^2}{\theta},
\end{equation}
with 
\begin{equation}\label{defLambda}
\Lambda(\theta) := \int_1^\theta \frac{1}{s}\dd s = \log\theta.
\end{equation}
We remark that in our case $\Lambda(\theta)$ is a very well-known function, but we stick with this notation in order to be coherent with the literature \cite{m.eleuterie.roccag.schimperna2015, ms}, where $\Lambda(\theta)$ might be a more generic function.
Testing \eqref{Leq} by $\xi \in C^\infty([0,T]\times \overline{\Omega}),\ \xi\geq0,\ \xi( T,\cdot)=0$ and integrating by parts we infer
\begin{align*}
&\itt \io (\Lambda(\theta) + \fhi) \xi_t \dd x \dd t + \itt\io \dfrac{\kappa(\theta)}{\theta}\nabla \theta \cdot \nabla \xi \dd x \dd t\\ 
 & \quad = -\itt \io \dfrac{|\nabla\mu|^2}{\theta} \xi \dd x \dd t -\itt \io \dfrac{\kappa(\theta)}{\theta^2}|\nabla \theta|^2 \dd x \dd t -\io (\Lambda(\theta_0) + \fhi_0) \xi(\cdot,0) \dd x.
\end{align*}
Setting $\delta(\theta): = \int_1^\theta \frac{\kappa(s)}{s} \dd s = \ln \theta +\frac{1}{q}(\theta^q -1)$ according to \eqref{heatassumption}, we get
\begin{align}\nonumber
&\itt \io (\Lambda(\theta) + \fhi) \xi_t \dd x \dd t + \itt\io \delta(\theta) \Delta \xi \dd x \dd t\\ \label{Lineq}
 & \quad = -\itt \io \dfrac{|\nabla\mu|^2}{\theta} \xi \dd x \dd t -\itt \io \dfrac{\kappa(\theta)}{\theta^2}|\nabla \theta|^2 \dd x \dd t -\io (\Lambda(\theta_0) + \fhi_0) \xi(\cdot,0) \dd x.
\end{align}

\begin{defn}\label{defensol}
We say that $(\fhi, \mu, \theta, \sigma)$ is a \textit{weak entropy solution} to our non-isothermal Cahn-Hilliard model if it sastisfies the following equations
\begin{align*}
&\langle \fhi_t,\xi\rangle = -\io \nabla \mu \cdot \nabla \xi \dd x +\io (\PP \sigma -\AA)h(\fhi)\xi \dd x\qquad \textrm{a.e. in }(0,T) \textrm{ and } \forall \xi \in V, \smallskip\\
&\mu = -\Delta \fhi + F'(\fhi) -\theta \qquad \textrm{a.e. in }(0,T)\times\Omega, \smallskip\\
&\langle \sigma_t,\xi\rangle = -\io \nabla \sigma \cdot \nabla \xi \dd x -\io \CC \sigma h(\fhi)\xi \dd x+\io \BB(\sigma_B-\sigma)\xi \dd x\\
&\textrm{a.e. in }(0,T) \textrm{ and } \forall \xi \in V,
\end{align*}
complying a.e. in $\Omega$ with the initial conditions \eqref{ic}, homogeneus Neumann boundary conditions and the entropy production inequality
\begin{small}
\begin{align} \nonumber
&\itt \io (\Lambda(\theta) + \fhi) \xi_t \dd x \dd t + \itt\io \delta(\theta) \Delta \xi \dd x \dd t\\ 
& \quad \leq -\itt \io \dfrac{|\nabla\mu|^2}{\theta} \xi \dd x \dd t -\itt \io \dfrac{\kappa(\theta)}{\theta^2}|\nabla \theta|^2 \dd x \dd t-\io (\Lambda(\theta_0) + \fhi_0) \xi(\cdot,0) \dd x 
\end{align}
\end{small}
$\forall \xi \in C^\infty([0,T]\times \overline{\Omega}),\ \xi\geq0,\ \xi( T,\cdot)=0$.
\end{defn}

\begin{thm}
\label{maint}
Suppose that the assumptions in Section \ref{ass} hold and let $T>0$. Then there exists at least one weak solution to our model problem, namely a quadruple $(\fhi, \mu, \theta, \sigma)$ with regularity
\begin{align*}
& \fhi \in C([0,T];V)\cap H^1(0,T;V')\cap L^2(0,T;H^2(\Omega))\\
& \beta(\fhi)\in L^2(0,T;H)\\
&\mu \in L^2(0,T;V)\\
&\theta \in L^2(0,T;V) \cap L^\infty(0,T;L^1(\Omega)) \cap L^q(0,T;L^{3q}(\Omega)), \; q\geq 2, \quad \theta >0 \textrm{ a.e. in }\Omega\\
& \sigma \in C([0,T];H)\cap H^1(0,T;V') \cap L^\infty((0,T)\times\Omega)\cap L^2(0,T;V)
\end{align*}
satisfying system \eqref{eqfhi}--\eqref{eqnutr} in the sense of Definition \ref{defensol}.
\end{thm}

\section{A priori estimates}
\label{ape}
This section is devoted to gain the suitable regularity for the quadruple $(\fhi,\mu,\theta,\sigma)$ to prove the existence of solutions in Section \ref{wss}. These a priori bounds are obtained formally, working directly on our system \eqref{eqfhi}-\eqref{eqnutr}. We remark that the existence (of weak entropy solutions) argument might be made rigorous by the Faedo-Galerkin method that we decided not to detail here.

\subsection{Nutrient estimate}
We first search for a priori bounds for the nutrient following \cite{mrs19}. Therefore we give here only a sketch of the main steps.\\
Testing \eqref{eqnutr} by $-\sigma_-$, where $\sigma_- \geq 0$ represents the negative part of the nutrient $\sigma$, exploting the initial conditions on $\sigma$ and applying the Gronwall lemma, we gain
\[ \sigma(t,x) \geq 0 \quad \textrm{for a.e. } t\geq 0, x\in \Omega.\] 
Now, testing \eqref{eqnutr} by $(\sigma- \bar{\sigma})_+$ (where $\bar{\sigma}\geq 1$ is a suitable constant) using the Gronwall lemma and our assuptions on $h$ and $\sigma_B$, it is possible to obtain
\begin{align}\label{nutrest}
\| \sigma \|_{L^\infty ((0,T)\times \Omega)}\leq c_T,
\end{align}
where $c_T$ is a constant depending on time.

\subsection{Energy estimate}
We test \eqref{eqfhi} by $\mu$, \eqref{eqmu} by $\fhi_t$ and \eqref{eqtemp} by 1 and then sum up. This yields, taking into account the boundary conditions,
\begin{align}\label{en1}
\ddt \left(\12\Vert \nabla \fhi\Vert^2 +\io F(\fhi) +\io \theta\right) = \io(\PP\sigma-\AA)h(\fhi)\mu.
\end{align}

We take care of the right hand side, in particular
\begin{small}
\begin{align*}
&\io(\PP\sigma-\AA)h(\fhi)\mu \\
\stackrel{\eqref{eqmu}}{=} &-\io (\PP\sigma-\AA) h(\fhi)\Delta\fhi +\io(\PP\sigma-\AA)h(\fhi)F'(\fhi)- \io(\PP\sigma-\AA)h(\fhi)\theta\\
\stackrel{\eqref{F}}{=}&\io(\PP\sigma-\AA)h'(\fhi)|\nabla\fhi|^2+\io \PP h(\fhi)\nabla\sigma\cdot\nabla\fhi+\io\beta(\fhi)(\PP\sigma-\AA)h(\fhi)\\
&+ \io \lambda\fhi(\AA-\PP\sigma)h(\fhi)+\io(\AA-\PP\sigma)h(\fhi)\theta.
\end{align*}
\end{small}
Thus \eqref{en1} reads
\begin{small}
\begin{align} \nonumber
\ddt \left(\12\Vert \nabla \fhi\Vert^2 +\io F(\fhi) +\io \theta\right) = &\io(\PP\sigma-\AA)h'(\fhi)|\nabla\fhi|^2\\ \nonumber
&+\io(\PP\sigma-\AA)\beta(\fhi)h(\fhi)\\ \nonumber
&+ \io \lambda(\AA-\PP\sigma)h(\fhi)\fhi+\PP\io h(\fhi)\nabla\sigma\cdot\nabla\fhi\\ \nonumber
&+ \io(\AA-\PP\sigma)h(\fhi)\theta\\\label{en2}
:= & I+II +III+IV.
\end{align}
\end{small}
We now estimate each term on the right hand side separately. The estimate on the nutrient \eqref{nutrest} is a key point for all these bounds. In particular this is where a time-dependent constant $c_T$ comes from. Exploting the assumption \eqref{hh'}, we infer
\begin{align}\label{estI}
I \leq c_T\|\nabla \fhi\|^2.
\end{align} 
According to \eqref{beta} it is straightforward that
\begin{align}\label{estII}
II \leq c_T\left(1+\io F(\fhi)\right).
\end{align}
Moreover, using once again the assumption \eqref{hh'} on $h$ and Young's inequality, we get
\begin{align}\label{estIII}
III\leq \12 \|\nabla \sigma\|^2 + c_T \left(1+\|\fhi\|_{L^1(\Omega)} + \|\nabla \fhi\|^2 \right).
\end{align}
Eventually, by the same tools used to estimate $III$, it holds
\begin{align}\label{estIV}
IV \leq c_T \|\theta \|_{L^1(\Omega)}.
\end{align}
Combining estimates \eqref{estI}--\eqref{estIV}, \eqref{en2} reads
\begin{align}\label{en3}
&\ddt \left(\12\Vert \nabla \fhi\Vert^2 +\io F(\fhi) +\io \theta\right)\\ \nonumber
\leq \ & \12 \|\nabla \sigma\|^2 + c_T \left(1+\|\fhi\|_{L^1(\Omega)} + \|\nabla \fhi\|^2 + \io F(\fhi)+ \|\theta \|_{L^1(\Omega)}\right)
\end{align}
Our aim is to apply Gronwall's lemma in order to gain the energy estimate. Therefore we estimate and reabsorb the term $\|\fhi\|_{L^1(\Omega)}$ according to \eqref{bsup}. Moreover we test \eqref{eqnutr} by $\sigma$ which yields
\begin{align}
\12 \ddt \|\sigma\|^2 +\|\nabla \sigma\|^2\leq c(1+\|\sigma\|^2).
\end{align}
Hence, summing this last estimate to \eqref{en3} we finally get
\begin{align}\nonumber
&\ddt \left(\12\Vert \nabla \fhi\Vert^2 +\io F(\fhi) +\io \theta + \12 \|\sigma\|^2\right)+ \12 \|\nabla \sigma\|^2\\\label{en4}
\leq\ & c_T \left(1+\|\nabla \fhi\|^2 + \io F(\fhi)+ \|\theta \|_{L^1(\Omega)} + \|\sigma\|^2\right).
\end{align}
We are now able to apply Gronwall's lemma to \eqref{en3}, therefore we obtain the following a priori estimates
\begin{align}\label{Ggrfhi}
& \|\nabla \fhi\|_{L^\infty(0,T;H)}\leq c_T\\ \label{GF}
& \|F(\fhi)\|_{L^\infty(0,T;L^1(\Omega))} \leq c_T\\ \label{Gtemp}
&\|\theta\|_{L^\infty(0,T;L^1(\Omega))}\leq c_T \\ \label{Gnutr}
&\|\sigma\|_{L^\infty(0,T;H)\cap L^2(0,T;V)}\leq c_T .
\end{align}

\noindent In particular, combining \eqref{beta} and \eqref{bsup} with \eqref{GF}, we gain
\begin{equation}\label{entrfhi}
\|\fhi\|_{L^\infty(0,T;L^1(\Omega))} \leq c_T
\end{equation}
\subsection{Entropy estimate}
We now derive the entropy estimate testing \eqref{eqtemp} by $-\dfrac{1}{\theta}$. Therefore
\begin{equation}\label{entrest}
  \ddt \io ( - \log \theta - \fhi )
   + \io \dfrac{1}{\theta} | \nabla \mu |^2
   + \io \big( | \nabla \log \theta |^2 + k_q | \nabla \theta^{q/2} |^2 \big)
   = 0,
\end{equation}

where $k_q>0$ is a suitable constant only depending on the exponent $q\geq 2$, introduced in \eqref{heatassumption}.\\Now, integrating in time, owing to \eqref{Gtemp} and \eqref{entrfhi} and recalling that $|\log r| \leq r - \log r \;\; \forall r > 0$, we infer
\begin{align}\label{entrlog}
  & \| \log \theta \|_{L^\infty(0,T;L^1(\Omega))} 
    + \| \log \theta \|_{L^2(0,T;V)} \le c_T,\\ \label{entrgradtemp}
  & \| \nabla \theta^{q/2} \|_{L^2(0,T;H)} \le c_T.
\end{align}  
Then, combining \eqref{nlp} with \eqref{Gtemp} and \eqref{entrgradtemp}, it holds
\begin{align}\label{tempv}
\| \theta^{\frac{q}{2}}\|_{L^2(0,T;V)}\leq c_T
\end{align}
which implies in particular, since $q\geq 2$
\begin{align}\label{tempest}
\| \theta\|_{L^2(0,T;V)}\leq c_T.
\end{align}
On the other hand, using Sobolev embedding theorems, \eqref{tempv} also implies
\begin{align*}
\| \theta^{\frac{q}{2}}\|_{L^2(0,T;L^6(\Omega))}\leq c_T
\end{align*}
and hence
\begin{align}\label{tempq}
\| \theta\|_{L^q(0,T;L^{3q}(\Omega))}\leq c_T.
\end{align}

\subsection{Chemical potential estimate}
Integrating \eqref{eqtemp} over $\Omega$ and exploiting boundary conditions together with Gauss-Green formula, we infer
\begin{align}\label{pottemp}
\|\nabla \mu\|^2 = \ddt \io \theta +\io \theta \fhi_t.
\end{align}
We now rewrite the latter term according to \eqref{eqfhi}, then using \eqref{hh'} and \eqref{nutrest}, it follows that \eqref{pottemp} reads
\begin{align}
\12 \| \nabla \mu\|^2 \leq \ddt \io \theta + \12 \| \nabla \theta\|^2 +c_T \|\theta\|_{L^1(\Omega)}.
\end{align}
Thus from \eqref{Gtemp} and \eqref{tempest}, we obtain
\begin{align}\label{gradpot}
\|\nabla \mu\|_{L^2(0,T;H)}\leq c_T.
\end{align}
Now we integrate \eqref{eqmu} over $\Omega$, then 
\begin{align}
|\mu_\Omega| &\stackrel{\eqref{mean},\eqref{F}}{=} \Big\vert \io (\beta(\fhi)-\lambda \fhi) -\io \theta\Big\vert \ \\
& \quad \;\leq\io |\beta(\fhi)|+\io |\lambda\fhi|+\|\theta\|_{L^1(\Omega)}\\
&\stackrel{\eqref{beta},\eqref{Gtemp},\eqref{entrfhi}}{\leq} c_{\beta}\left(1+\io F(\fhi)\right) + c_T.
\end{align}
Using now the bound \eqref{GF}, we get
\begin{align}
\|\mu_\Omega\|_{L^\infty(0,T)} \leq c_T.
\end{align}
Combining this last bound with the Poincar\'{e} inequality \eqref{poin} and the previous estimate \eqref{gradpot}, we achieve 
\begin{align}\label{estpot}
\| \mu\|_{L^2(0,T;V)}\leq c_T.
\end{align}

\subsection{$\fhi$-dependent estimates}
We start testing \eqref{eqfhi} by $\fhi$, which leads to
\begin{align}
\ddt \io |\fhi|^2 = -\io \nabla \mu \cdot \nabla \fhi +\io (\PP \sigma -\AA)h(\fhi)\fhi.
\end{align}
Exploiting Young's inequality, the uniform bounds on $h$ and \eqref{nutrest} we infer
\begin{align*}
\ddt \|\fhi\|^2 \leq \12 \| \nabla \mu\|^2 +\12 \|\nabla \fhi\|^2 + c_T \|\fhi\|_{L^1(\Omega)}.
\end{align*}
Thus, integrating in time and using \eqref{gradpot}, \eqref{Ggrfhi} and \eqref{entrfhi} we get
\begin{align*}
\|\fhi\|_{L^\infty(0,T;H)}\leq c_T,
\end{align*}
whence estimate \eqref{Ggrfhi} gives
\begin{align}\label{estfhi}
\|\fhi\|_{L^\infty(0,T;V)}\leq c_T.
\end{align}
Next we test \eqref{eqmu} by $\beta(\fhi)$ and we obtain
\begin{align*}
\io |\beta(\fhi)|^2 + \io \beta'(\fhi)|\nabla \fhi|^2 = \io \mu \beta(\fhi) + \io \lambda \fhi \beta(\fhi) + \io \theta \beta(\fhi)
\end{align*}
Now, from \eqref{estpot}, \eqref{estfhi}, \eqref{tempest} and the monotonicity of $\beta$, it follows
\begin{align}\label{estbeta}
\|\beta(\fhi)\|_{L^2(0,T;H)}\leq c_T.
\end{align}
Taking advantage of this last estimate with \eqref{beta} and again of \eqref{estfhi} and \eqref{tempest}, a direct comparison within equation \eqref{eqmu} yields
\begin{align}\label{estfhiell}
\|\fhi\|_{L^2(0,T;H^2)}\leq c_T
\end{align}

\subsection{Further regularity}
We start testing \eqref{eqfhi} by a nonzero test function $v \in V$ and we infer
\begin{align*}
\langle\fhi_t,v\rangle = - \io \nabla\mu \cdot \nabla v + \io (\PP\sigma-\AA)h(\fhi)v.
\end{align*}
Now, according to estimates \eqref{nutrest}, \eqref{gradpot} and \eqref{estfhi} it follows
\begin{align}\label{estfhit}
\|\fhi_t\|_{L^2(0,T;V')}\leq c_T.
\end{align}
Taking advantage of this last estimate and exploting \eqref{estfhiell} together with \eqref{estfhi}, we infer (for example from \cite{lionsmag})
\begin{equation}
\fhi \in C([0,T];V).
\end{equation}

Similarly, multiplying equation \eqref{eqnutr} by a nonzero test function $v \in V$ and exploiting the bound \eqref{Gnutr}, it holds
\begin{align}\label{estnutrt}
\|\sigma_t\|_{L^2(0,T;V')}\leq c_T.
\end{align}
From standard embedding results (see e.g. \cite{brezzigilardi}), combining \eqref{estnutrt} and \eqref{Gnutr}, we gain the additional regularity for the nutrient
\begin{equation}
\sigma \in C([0,T];H).
\end{equation}

\section{Weak sequential stability}
\label{wss}
We assume to have a sequence of weak solutions $(\fhi_n, \mu_n, \theta_n,\sigma_n)$ which satisfies the a priori estimates obtained in Section \ref{ape} uniformly with respect to $n \in \mathbb{N}$. \\
We then show that, by weak compactness arguments, up to the extraction of a subsequence, $(\fhi_n, \mu_n, \theta_n,\sigma_n)$ converges in a suitable way to an entropy solution to our problem, i.e., to a limit quadruple $(\fhi, \mu, \theta,\sigma)$ solving \eqref{eqfhi}--\eqref{eqnutr} in the sense of Theorem \ref{maint}.\\
Indeed, exploiting the above bounds \eqref{nutrest}, \eqref{Gtemp}, \eqref{Gnutr}, \eqref{tempest}, \eqref{tempq}, \eqref{estpot}, \eqref{estfhi}, \eqref{estfhiell}, \eqref{estfhit} and \eqref{estnutrt}, together with standard weak compactness results, it is possible to extract a nonrelabelled subsequence such that
\begin{footnotesize}
\begin{align}\label{wfhi}
& \fhi_n \rightarrow \fhi \textrm{ weakly star in } L^\infty(0,T;V)\cap L^2(0,T;H^2(\Omega))\cap H^1(0,T;V')\\\label{wmu}
&\mu_n \rightarrow \mu  \textrm{ weakly in }L^2(0,T;V)\\ \label{wtemp}
&\theta_n \rightarrow \theta \textrm{ weakly star in } L^2(0,T;V) \cap L^\infty(0,T;L^1(\Omega)) \cap L^q(0,T;L^{3q}(\Omega))\\ \label{wnutr}
&\sigma_n \rightarrow \sigma \textrm{ weakly star in } L^\infty(0,T)\times\Omega))\cap L^2(0,T;V)\cap H^1(0,T;V')
\end{align}
\end{footnotesize}
Moreover combining  \eqref{estfhit} and \eqref{estnutrt} with \eqref{wfhi} and \eqref{wnutr} respectively and applying the Aubin-Lions lemma, we infer that
\begin{eqnarray}\label{convfhinutr}
\fhi_n \rightarrow \fhi \textrm{ and } \sigma_n \rightarrow \sigma \textrm{ strongly in } L^2(0,T;H).
\end{eqnarray}
Moreover convergence \eqref{wtemp} and interpolation theory for $L^p$ spaces imply that 
\begin{equation}\label{convtemp}
\theta_n \rightarrow \theta \textrm{ strongly in } L^p(0,T;L^p(\Omega)),\quad  p\in \left[1,q+\frac{2}{3}\right).
\end{equation} 

Indeed, from standard interpolation theory, we know that if $f\in L^p \cap L^s$, then $f \in L^r$, with $r$ s.t. $\dfrac{1}{r} = \dfrac{\gamma}{p}+ \dfrac{1-\gamma}{s}$. We first consider the time-spaces $L^\infty$ and $L^q$, hence $s = \infty$ and $p=q$ which gives
\begin{equation}\label{interp1}
\frac{1}{r} = \dfrac{\gamma}{q}.
\end{equation}
We then apply the general interpolation result to the space-spaces $L^1$ and $L^{3q}$, from which it follows
\begin{equation}\label{interp2}
\frac{1}{r} = \dfrac{\gamma}{3q}+\frac{1-\gamma}{1}.
\end{equation}
Since \eqref{interp1} and \eqref{interp2} must hold simultaneously, we infer that $r = q+\dfrac{2}{3}$.

Now, according to (Theorem 2.19, \cite{brezzigilardi}) with $s=r=0$ it follows that
\[L^\infty(0,T;L^1(\Omega)) \cap L^q(0,T;L^{3q}(\Omega))\hookrightarrow\hookrightarrow L^p(0,T;L^p(\Omega))\quad  p\in \left[1,q+\frac{2}{3}\right).\]

Therefore it is possible to pass to the limit also in the nonlinear terms, according to the continuity of $\kappa, \beta$ and $h$. Indeed, by a generalized version of Lebesgue's dominated convergence theorem it holds
\begin{align}\label{convk}
&\kappa(\theta_n)\rightarrow \kappa(\theta) \textrm{ strongly in } L^p(0,T;L^p(\Omega)),\quad p\in \Big[1,1+\frac{2}{3q}\Big)\\
&\beta(\fhi_n)\rightarrow \beta(\fhi) \textrm{ weakly in } L^2(0,T;H).
\end{align}

We now want to pass to the limit in the balance of entropy. Namely let us assume that \eqref{Leq} is satisfied by the approximate solution $(\fhi_n,\mu_n,\theta_n,\sigma_n), \forall n\in\mathbb{N}$. Testing it by $\xi \in C^\infty([0,T]\times \overline{\Omega}),\ \xi\geq0, \ \xi( T,\cdot)=0$ and integrating by parts we infer
\begin{small}
\begin{align}\nonumber
&\itt \io (\Lambda(\theta_n) + \fhi_n) \xi_t \dd x \dd t + \itt\io \delta(\theta_n) \Delta \xi \dd x \dd t\\ \label{Leqn}
 & \quad = -\itt \io \dfrac{|\nabla\mu_n|^2}{\theta_n} \xi \dd x \dd t -\itt \io \dfrac{\kappa(\theta_n)}{\theta_n^2}|\nabla \theta_n|^2 \dd x\dd t -\io (\Lambda(\theta_0) + \fhi_0) \xi(\cdot,0) \dd x,
\end{align}
\end{small}
where $\delta(\theta_n): = \int_1^{\theta_n} \frac{\kappa(s)}{s} \dd s$.\\
We first take care of the terms on the left hand side. 
According to \eqref{defLambda}, by \eqref{entrlog} and \eqref{convtemp},
\begin{equation}\label{convL}
\Lambda(\theta_n) \rightarrow \Lambda(\theta_n) \quad \textrm{strongly in } L^{1+}(0,T)\times \Omega).
\end{equation}

Moreover, from \eqref{entrlog} and \eqref{convtemp}, it follows that
\begin{equation}\label{convh}
\delta(\theta_n)\rightarrow \delta(\theta) \quad \textrm{strongly in } L^{1+}((0,T)\times \Omega),
\end{equation}
hence in particular
\[\itt \io \delta(\theta_n)\Delta \xi \to \itt \io \delta(\theta)\Delta \xi. \]

Then the first row of \eqref{Leqn} passes to the desired limit not only as a supremum limit, but as a true limit.
In order to deal with the first two terms in the right hand side we recall a useful lower semicontinuity result by Ioffe.
\begin{thm}[Ioffe]
Let $\mathcal{O}\subset \real^d$ be a smooth bounded open set and $f:\mathcal{O}\times \real^n \times \real^m \rightarrow
[0,+\infty], d, n, m \geq 1,$ be a measurable non-negative function such that
\begin{itemize}
\item[] $f(x,\cdot, \cdot)$ is lower semicontinuous on $\real^n\times \real^m$ for every $x\in \mathcal{O}$,
\item[]$f(x, u, \cdot)$ is convex on $\real^m$ for every $(x,u)\in \mathcal{O} \times \real^n$.
\end{itemize}
Let also $(u_k, v_k), (u,v): \mathcal{O} \rightarrow \real^n\times \real^m$ be measurable functions s. t. \smallskip\\
$u_k(x) \to u(x)$ in measure in $\mathcal{O},\qquad v_k \rightharpoonup v$ weakly in $L^1(\mathcal{O};\real^m)$.\smallskip\\
Then,
\[\liminf_{k\to +\infty} \int_{\mathcal{O}}
f(x, u_k(x), v_k(x)) \geq
 \int_{\mathcal{O}}f(x, u(x), v(x)).\]
\end{thm}
We start considering the first term in the right hand side. We exploit this result setting $\mathcal{O} =\Omega \times (0,T), f:\mathcal{O}\times \real^+ \times \real ^3 \rightarrow [0,\infty] \textrm{ s.t } (x,t)\times w \times v \mapsto w|v|^2$. Such $f$ satisfies Ioffe's assumptions.
 putting $w_n = \frac{\xi}{\theta},\  v_n = \nabla \mu_n, \ \forall n\in \mathbb{N}$. Hence, by \eqref{wmu} it holds $\lbrace \nabla \mu_n\rbrace_n \rightharpoonup \nabla \mu$ in $L^1(\mathcal{O})$. Therefore by Ioffe's theorem,
\begin{equation}\label{ioffemu}
 \liminf_{n\to +\infty} \itt \io |\nabla \mu_n|^2\xi \geq \itt \io |\nabla \mu|^2\xi.
 \end{equation}
In a similar way, from \eqref{convtemp} and \eqref{convk},
\begin{equation}\label{ioffek}
\liminf_{n\to +\infty} \itt \io \dfrac{\xi}{\theta_n}\dfrac{\kappa(\theta_n)}{\theta_n}|\nabla \theta_n|^2 \geq  \itt \io \dfrac{\xi}{\theta}\dfrac{\kappa(\theta)}{\theta}|\nabla \theta|^2.
\end{equation}

Furthermore, assuming that $\theta_n(0,\cdot)$ converges in a suitable way to $\theta_0$, putting together \eqref{convfhinutr}, \eqref{convL}, \eqref{convh}, \eqref{ioffemu} and \eqref{ioffek}, it follows that we eventually recover \eqref{Lineq}. It is worth noting that the inequality sign is due to the application of Ioffe's theorem.
This concludes the procedure and so the proof of existence of weak entropy solutions.

\begin{rmk}
We notice that we have assumed throughout the proof that the absolute temperature is a.e. positive. This is crucial in order for estimates in Section \ref{ape} to make sense. In particular it should be shown that the solution $\theta_n$ of the discretized problem (for instance in a Faedo-Galerkin scheme, that we decided not to detail here) is positive.
At least, according to \eqref{entrlog} the strict positivity of $\theta_n$ will be preserved a.e. in $(0,T)\times\Omega$ also in the limit.
\end{rmk}

\medskip
\smallskip

\noindent\textbf{Acknowledgements}\\
The author is supported by GNAMPA (Gruppo Nazionale per l'Analisi Matematica, la Probabilit\`{a} e le loro Applicazioni) of INdAM (Istituto Nazionale di Alta Matematica) and by MIUR through the project FFABR (M. Eleuteri). The author would also like to thank Prof. Michela Eleuteri for careful reading and helpful comments.

\end{document}